\newcommand{\R}{\mathbb{R}}
\newcommand{\ds}{\displaystyle} 
\mathchardef\varepsilon="010F
\mathchardef\epsilon="0122
\mathchardef\vartheta="0112
\mathchardef\theta="0123
\mathchardef\varrho="011A
\mathchardef\rho="0125
\mathchardef\varphi="011E     
\mathchardef\phi="0127
\renewcommand \emptyset \varnothing
\documentclass[11pt]{article}            
\usepackage[all]{xy}
\usepackage[latin1]{inputenc}        
\usepackage[dvips]{graphics}
\usepackage[dvips]{graphicx}
\usepackage{amsfonts}
\usepackage{amssymb}
\usepackage{amsmath}
\usepackage{amstext}
\usepackage{amsbsy}
\usepackage{amsopn}
\usepackage{amsthm}
\usepackage{amscd}
\usepackage{amsxtra}
\usepackage{mathrsfs}
\usepackage{upref}
\date{}                            
\title{An Obstacle Control problem involving the $p$-Laplacian}

\author{\sc Mouna Kraïem\\
{\small Université de Cergy Pontoise , 
Departement de Mathematiques, Site de Saint-Martin} \\
{\small 2, avenue Adolphe Chauvin, 95302 Cergy Pontoise Cedex, France} \\ 
{\small e-mail: Mouna.Kraiem@math.u-cergy.fr}}

\numberwithin{equation}{section}                
\theoremstyle{plain}
\newtheorem{theorem*}{\bf Theorem}
\newtheorem{theorem}{\bf Theorem}
\newtheorem{lemma}{\bf Lemma}
\newtheorem{proposition}{\bf Proposition}
\newtheorem*{proposition*}{\bf Proposition}

\newtheorem*{corollary*}{\bf Corollary}
\theoremstyle{definition}
\newtheorem*{definition*}{\bf Definition}
\newtheorem{definition}{\bf Definition}

\newtheorem*{example*}{\bf Example}
\theoremstyle{remark}

\begin{document}

\maketitle

\begin{abstract}
In this paper, we consider the analogous of the obtacle problem in 
$H_0^1(\Omega)$, on the space $W^{1,p}_0(\Omega)$. 
We prove an existence and uniqueness of the result.
In a second time, we define the optimal control problem associated. 
The results, here enclosed, generalize the one obtained by D.R. Adams, S. Lenhard
in \cite{1}, \cite{2} in the case $p=2$.
\end{abstract}


\section{Introduction}
Let $\Omega$ be a bounded domain in $\R^N, N\geq2$, whose boundary is $\mathcal C^1$
piecewise. For $p >1$ and for $\psi$ given in 
$W^{1,p}_0(\Omega)$, define
$$K(\psi)=\{v \in W_0^{1,p}(\Omega), v \geq \psi \ \text{a.e. in} \ \Omega\}.$$ 
It is clear that $K(\psi)$ is a convex and weakly closed set in $L^{p}(\Omega)$. 
Let $p'$ be the conjugate of $p$, and $f \in L^{p'}(\Omega)$. 
We consider the following variational inequality called {\it the obstacle problem}:
\begin{equation}\label{P-ob:eq:1.1}
\begin{cases}
u \in K(\psi),\\
{\ds \int_\Omega} \sigma(u) \cdot \nabla(v-u) \ dx \geq {\ds \int_\Omega} f(v-u) \ dx, 
\ \forall v \in K(\psi),
\end{cases}
\end{equation}
where $\sigma(u)= |\nabla u|^{p-2} \nabla u$.
We shall say that $\psi$ is {\it the obstacle} and $f$ is {\it the
source term}.

We begin to prove existence and uniqueness of a solution $u$ to (\ref{P-ob:eq:1.1}), 
using variational formulation of the obstacle problem on the set $K(\psi)$.
We shall then denote $u$ by: $u=T_f(\psi)$.
Secondly, we characterize $T_f(\psi)$ as the lowest $f-$superharmonic
function greater than $\psi$.

\section{Existence and  uniqueness of the solution}
\begin{proposition}\label{P-ob:prop:1.1}
A function $u$ is a solution to the problem {\rm{(\ref{P-ob:eq:1.1})}} if and only if $u$
satisfies the following:
\begin{equation}\label{P-ob:eq:1.2}
\begin{cases}
u \in K(\psi),\\
-\Delta_p u \geq f, \  \text{a.e. \ in} \ \Omega,\\
{\ds \int_\Omega }\sigma(u) \cdot \nabla (\psi-u) \ dx={\ds \int_\Omega} f(\psi-u) \ dx.
\end{cases}
\end{equation}
\end{proposition}

\begin{proof}[Proof of Proposition \ref{P-ob:prop:1.1}]
Suppose that $u$ satisfies (\ref{P-ob:eq:1.1}). Then taking $v=u+\phi\in
K(\psi)$ for $\phi \in \mathcal D(\Omega), \ \phi \geq 0$, one gets that
$-\Delta_p u \geq f$ in $\Omega$.
Moreover, For $v=\psi$ and $v= 2u-\psi$, one gets that
$$\int_\Omega \sigma(u) \cdot \nabla (\psi-u) \ dx=\int_\Omega f(\psi-u) \ dx,$$
hence $u$ satisfies (\ref{P-ob:eq:1.2}).
 
Conversely, let $u \in K(\psi)$ such that $-\Delta_p u \geq
f$, let $v$ be in $K(\psi)$ and  $\phi_n \in
\mathcal D(\Omega), \phi_n \geq 0$ such that 
$\phi_n \rightarrow v-\psi \ \text{in} \ W^{1,p}_0(\Omega)$. Then one gets
\begin{align*}
\int_\Omega \sigma(u) \cdot \nabla (v-\psi)
&=\lim_{n\rightarrow \infty} \int_\Omega \sigma(u) \cdot \nabla \phi_n\\
&=\lim_{n\rightarrow \infty} \int_\Omega -\Delta_p u \  \phi_n\\
&\geq \lim_{n\rightarrow \infty} \int_\Omega f\phi_n =
\int_\Omega f (v-\psi), \ \forall \ v\in K(\psi).
\end{align*}
Using the last equality of (\ref{P-ob:eq:1.2}), one gets that 
$$\int_\Omega \sigma(u) \cdot \nabla (v-u) \geq \int_\Omega f (v-u), 
\ \forall \ v\in K(\psi),$$
hence $u$ satisfies (\ref{P-ob:eq:1.1}).
\end{proof}

Let us prove now the existence and uniqueness of a solution to the obstacle problem (\ref{P-ob:eq:1.1}).

\begin{proposition}\label{P-ob:prop:1.2}
There exists a solution to {\rm{(\ref{P-ob:eq:1.1})}}, which can be obtained as the 
minimizer of the following minimization problem
\begin{equation}\label{P-ob:eq:1.3}
\inf_{v \in K(\psi)} I(v),
\end{equation}
where $I$ is the following energy functional
$$I(v)=\frac{1}{p}\int_\Omega |\nabla v|^p -
\int_\Omega f v.$$
\end{proposition}

\begin{proof}[Proof of Proposition \ref{P-ob:prop:1.2}]  
Using classical arguments in the calculus of variations, since $K(\psi)$
is a weakly closed convex set in $W^{1,p}_0(\Omega)$, 
and the functional $I$ is convex and coercive on $W^{1,p}_0(\Omega)$, then 
one obtains that there exists a solution $u$ to (\ref{P-ob:eq:1.3}). 
\end{proof}

\begin{proposition}\label{P-ob:prop:1.3}
The inequation {\rm{(\ref{P-ob:eq:1.1})}} possesses a unique solution.
\end{proposition}

\begin{proof}[Proof of Proposition \ref{P-ob:prop:1.3}]
Suppose that $u_1, u_2 \in W^{1,p}_0(\Omega)$ are two solutions of the variational 
inequality (\ref{P-ob:eq:1.1})
$$u_i \in K(\psi): 
\int_\Omega \sigma(u_i) \cdot \nabla(v-u_i) \ dx \geq \int_\Omega
f(v-u_i) \ dx, \ \forall \ v \in K(\psi), \
i=1,2$$
Taking $v=u_1$ for $i=2$ and $v=u_2$ for $i=1$ and adding, we have
$$\int_\Omega \left[\sigma(u_1)-\sigma(u_2)\right] \cdot  \nabla(u_1-u_2) \leq 0.$$
Recall that we have
$$\int_\Omega \left[\sigma(u_1)-\sigma(u_2)\right] \cdot  \nabla(u_1-u_2)\geq 0,$$
which implies that 
$$\int_\Omega \left[\sigma(u_1)-\sigma(u_2)\right] \cdot  \nabla(u_1-u_2)=0,$$
and then, $u_1=u_2$ a.e in $\Omega$.
\end{proof}

Thus, we get the existence and uniqueness of a solution to (\ref{P-ob:eq:1.1}).

\begin{definition}\label{P-ob:def:1.1}
We shall say that $u$ is {\it $f-$superhamonic} in $\Omega$,
if $u \in W^{1,p}_0(\Omega)$ is a weak solution to 
$-\Delta_p u \geq f$, in the sense of distributions.
\end{definition}

\begin{proposition}\label{P-ob:prop:1.4}
A function $u$ is a solution of  {\rm{(\ref{P-ob:eq:1.1})}}, if and only if $u$ is the lowest
$f-$superharmonic function, greater than $\psi$.
\end{proposition}

\begin{proof}[Proof of Proposition \ref{P-ob:prop:1.4}]
Let $u$ be a solution of (\ref{P-ob:eq:1.1}) and $v$ be an $f-$super\\harmonic
function, greater than $\psi$.
Let $\xi=\max(u,v), \ \xi \in K(\psi)$.
Recalling that $v^-=\sup(0, -v)$, one has then $(\xi-u)=-(v-u)^-$.
From (\ref{P-ob:eq:1.1}), one gets
$$\int_\Omega \sigma(u) \cdot \nabla (\xi-u) \geq \int_\Omega f
(\xi-u).$$
On the other hand, since $\xi-u \leq 0$ and $-\Delta_p v \geq f$, we have 
$$\int_\Omega \sigma(v) \cdot \nabla (\xi-u) \leq \int_\Omega f (\xi-u).$$
We obtain, subtracting the  above two inequalities:
$$\int_\Omega \left[\sigma(v)-\sigma(u)\right] \cdot \nabla (\xi-u) \leq 0,$$
which implies that
$$-\int_\Omega \left[\sigma(v)-\sigma(u)\right] \cdot \nabla (v-u)^- \leq 0,$$
and then $(v-u)^-=0$, or equivalently $u\leq v$ in $\Omega$.
\end{proof}
Recall that we define by $T_f(\psi)$ the lowest
$f-$superharmonic function, greater than $\psi$.
\begin{lemma}\label{P-ob:lemma:1.1}
The mapping $\psi \mapsto T_f(\psi)$ is increasing.
\end{lemma}

\begin{proof}[Proof of Lemma \ref{P-ob:lemma:1.1}]
Let $u_1=T_f(\psi_1)$ and $u_2=T_f(\psi_2)$, 
which are respectively solutions to the following variational inequalities
$$\begin{cases}
-\Delta_p u_i \geq f\\
u_i \geq \psi_i, \ i=1,2
\end{cases}$$
and let $\psi_1 \leq \psi_2$. It is clear that $u_2 \geq \psi_1$.
Hence $u_2$ is $f-$superharmonic and using Proposition \ref{P-ob:prop:1.4}, 
one obtains $u_1\leq u_2$. 
\end{proof}

\begin{proposition}\label{P-ob:prop:1.5}
The mapping $\psi \mapsto T_f(\psi)$ is weak lower semicontinuous,
in the sense that:
\begin{itemize}
\item If $\psi_k \rightharpoonup \psi$ weakly in
$W^{1,p}_0(\Omega)$, then 
$T_f(\psi) \leq {\ds \liminf_{k \rightarrow \infty}} T_f(\psi_k)$.
\item ${\ds \int_\Omega} |\nabla(T_f(\psi))|^p \leq {\ds \liminf_{k \rightarrow \infty}} {\ds \int_\Omega} |\nabla(T_f(\psi_k))|^p$.
\end{itemize}
\end{proposition}

\begin{proof}[Proof of Proposition \ref{P-ob:prop:1.5}]
Let $(\psi_k)$ be a sequence in $W^{1,p}_0(\Omega)$ which converges weakly
in $W^{1,p}_0(\Omega)$ to $\psi$, 
and let $\phi_k= {\rm{min}}(\psi_{k},\psi)$. 
Since $T_f$ is increasing, one gets that $T_f(\phi_k) \leq T_f(\psi_k)$. 
We now prove that $T_f(\phi_k)$ converges strongly in $W^{1,p}_0(\Omega)$
towards $T_f(\psi)$. This will imply that 
 $$T_f(\psi)= \lim _{k \rightarrow \infty}T_f(\phi_k) \leq 
\liminf_{k \rightarrow \infty} T_f(\psi_{k}).$$
We denote  $u_k$ as $T_f(\phi_k)$. It is clear that $u_k$ is
bounded in $W^{1,p}_0(\Omega)$ since $\phi_k \leq \psi$. 
Hence for a subsequence, still denoted
$u_k$, there exists some $u$ in $W^{1,p}_0(\Omega)$ such that
\begin{equation}\label{P-ob:eq:1.4}
\nabla u_k \rightharpoonup \nabla u \  {\rm{weakly \ in}} \
L^p(\Omega) , \ u_k \rightarrow u \ {\rm{strongly \ in}} \ L^p(\Omega).
\end{equation}
On the other hand, using the fact that $\phi_k$ 
converges weakly to $\psi$ in $W^{1,p}_0(\Omega)$ (see
Lemma \ref{P-ob:lemma:1.2} below), 
one gets the following assertion:
$$u_k \geq \phi_k \Longrightarrow u \geq \psi.$$
Let us prove now that $u$ is a solution of the minimizing problem (\ref{P-ob:eq:1.3}). 
For that aim, for $v\in K(\psi)$, since  $v \geq \psi \geq \phi_k$, we have
\begin{align*}
\frac{1}{p} \int_\Omega |\nabla u|^p -\int_\Omega f u
&\leq \liminf_{k \rightarrow \infty}\frac{1}{p} \int_\Omega 
|\nabla u_k|^p -\int_\Omega f u_k\\
&\leq \liminf_{k \rightarrow \infty}
 \inf_{w \geq \phi_k}\left\{\frac{1}{p} \int_\Omega
|\nabla w|^p -\int_\Omega f w\right\}\\
&\leq \frac{1}{p} \int_\Omega |\nabla v|^p -\int_\Omega f v.
\end{align*}
Then $u$ realizes the infimum in (\ref{P-ob:eq:1.3}). 
At the same time, since $u \in K(\psi)$, 
one has the following convergence
$$\frac{1}{p} \int_\Omega |\nabla u_k|^p -\int_\Omega f u_k \longrightarrow
\frac{1}{p} \int_\Omega |\nabla u|^p -\int_\Omega f u,
 \ {\rm{when}} \ k \rightarrow \infty,$$
which implies that $u_k$ converges strongly to $u$ in $
W^{1,p}_0(\Omega)$. We can conclude that $T_f(\phi_k)$ converges 
strongly to $T_f(\psi)$. 
\end{proof}

\begin{lemma}\label{P-ob:lemma:1.2}
Suppose that $\psi_k$ converges weakly to some $\psi$ in
$W^{1,p}_0(\Omega)$. Then, $\phi_k=\min(\psi_k,\psi)$ 
converges weakly to $\psi$ in $W^{1,p}_0(\Omega)$.
\end{lemma}

\begin{proof}[Proof of Lemma \ref{P-ob:lemma:1.2}]
We have 
$$\psi_k \longrightarrow \psi \quad \text{in} \ L^p(\Omega).$$
Then   
$$\phi_k=\frac{\psi_k+\psi-|\psi_k-\psi|}{2} \longrightarrow \psi \quad \text{in} \ L^p(\Omega).$$
Let us prove now that $|\nabla \phi_k|$ is bounded in 
$L^p(\Omega)$. For that aim, we write
\begin{align*}
\int_\Omega  |\nabla \phi_k|^p 
&=\int_\Omega \left|\nabla
\left(\frac{\psi_k+\psi-|\psi_k-\psi|}{2}\right)\right|^p\\
&\leq  C_p \left(\int_\Omega  |\nabla \psi_k|^p +\int_\Omega  |\nabla \psi|^p\right).
\end{align*}
Therefore the sequence $\phi_k$ is bounded in $W^{1,p}_0(\Omega)$, so it converges weakly, up to a subsequence, to $\psi$ in $W^{1,p}_0(\Omega)$. 
\end{proof}

\begin{proposition}\label{P-ob:prop:1.6}
The mapping $T_f$ is an involution, i.e. $T_f^2=T_f$.
\end{proposition}

\begin{proof}[Proof of Proposition \ref{P-ob:prop:1.6}]
Up to replacing $\psi$ by $u$ in the variational inequalities
(\ref{P-ob:eq:1.1}), and using proposition
\ref{P-ob:prop:1.4}, one gets that $u=T_f(u)$.
Then, we conclude that $ T_f^2(\psi)=T_f(\psi)$. 
\end{proof}

\section{A method of penalization}
Let $\mathcal M^+(\Omega)$ be the set of all nonnegative Radon
measures on $\Omega$ and $W^{-1,p^\prime}(\Omega)$ be the dual space of 
$W^{1,p}(\Omega)$ on $\Omega$ where $p^\prime$ is the conjugate of $p \
(1<p<\infty)$.
Suppose that $u$ solves (\ref{P-ob:eq:1.1}). Using the fact that a nonnegative 
distribution on $\Omega$ is a nonnegative measure on $\Omega$ (cf. \cite{DF9}), one gets the existence of 
$\mu \geq 0, \ \mu \in \mathcal M^+(\Omega)$, such that
\begin{equation}\label{P-ob:eq:1.5}
\int_\Omega \sigma(u) \cdot \nabla \Phi \ dx - \int_\Omega f \Phi \ dx=
\langle \mu,\Phi \rangle, \quad \forall \ \Phi \in \mathcal D(\Omega),
\end{equation}
that we shall also write 
$-\Delta_p u =f+\mu, \quad \mu\geq0 \  \text{in} \ \Omega.$

Let us introduce
\begin{equation}\label{P-ob:eq:1.6}
\beta(x)=
\begin{cases}
0, \quad  x >0,\\
x, \quad x\leq 0.
\end{cases}
\end{equation}
Clearly, $\beta$ is $\mathcal C^1$ piecewise, $\beta(x) \leq 0$ and is nondecreasing. 
Let us consider, for some $\delta >0$, the following semilinear elliptic equation:
\begin{equation}\label{P-ob:eq:1.7}
\begin{cases}
-\Delta_p u+ \frac{1}{\delta} \beta(u-\psi)=f,  \quad
\text{in} \ \Omega\\
u_{|{\partial \Omega}}=0.
\end{cases}
\end{equation}
We have the following existence result:

\begin{theorem}\label{P-ob:thm:2.1}
For any given $\psi \in W_0^{1,p}(\Omega)$ and $\delta>0$, {\rm(\ref{P-ob:eq:1.7})} possesses a unique
solution $u^\delta$. Moreover,
\begin{itemize}
\item[{\rm(1)}]$u^\delta \longrightarrow u \ {\rm{strongly \ in}} \ {W_0^{1,p}(\Omega)}$, as $\delta \longrightarrow 0$, with $u:=T_f(\psi)$.
\item[{\rm(2)}] There exists a unique $\mu \in W^{-1,p^\prime}(\Omega) \cap \mathcal M^+(\Omega)$ such that:
\begin{itemize}
\item[{\rm(i)}]$- \frac{1}{\delta} \beta(u^\delta-\psi) \rightharpoonup \mu \ {\rm{in}} \ W^{-1,p^\prime}(\Omega) \cap \mathcal M^+(\Omega).$
\item[{\rm(ii)}]$\langle \mu, T_f(\psi)-\psi \rangle=0$.
\end{itemize}
\end{itemize}
\end{theorem}

\begin{proof}[Proof of Theorem \ref{P-ob:thm:2.1}]
(1) Let $B$ be defined as $B(r)={\ds \int_0^r} \beta(s)  ds, \ \forall \ r\in \R$.
We introduce the following variational problem
\begin{equation}\label{P-ob:eq:1.8}
\inf_{v \in W_0^{1,p}(\Omega)} \left\{\frac{1}{p} \int_\Omega |\nabla v|^p +
\frac{1}{\delta} \int_\Omega B(v-\psi)-\int_\Omega f v \right\}.  
\end{equation}
The functional in (\ref{P-ob:eq:1.8}) is coercive, strictly convex and continuous.
As a consequence it possesses a unique solution
$u^\delta \in W_0^{1,p}(\Omega)$.
Since  $B(0)=0$, one has  
$$\frac{1}{p} \int_\Omega |\nabla u^\delta|^p+ \frac{1}{\delta} 
\int_\Omega B(u^\delta-\psi) - \int_\Omega f u^\delta 
\leq \frac{1}{p} \int_\Omega |\nabla \psi|^p - \int_\Omega f \psi,$$
since $B\geq 0$, then $u^\delta$ is bounded in $W_0^{1,p}(\Omega)$.
Extracting from $u^\delta$ a subsequence,
there exists $u$ in $W_0^{1,p}(\Omega)$, such that  
$$\nabla u^\delta \rightharpoonup \nabla u \  {\rm{weakly \ in}} \
L^p(\Omega), \ u^\delta \rightarrow u \ {\rm{strongly \ in}} \ L^p (\Omega).$$
Using $\frac{1}{\delta} {\ds \int_\Omega} B(u^\delta-\psi) \leq C$
and the continuity of $B$ one has
 $$0\leq  \int_\Omega B(u-\psi) \leq \liminf_{\delta \rightarrow 0}
\int_\Omega B(u^\delta-\psi)=0,$$ 
hence $u \in K(\psi)$.

We want to prove now that $u$ solves (\ref{P-ob:eq:1.1}).
Let $v \in K(\psi)$, since $B(r)\geq0, \ \forall \ r \in \R$ one gets:
\begin{align*}
\frac{1}{p} \int_\Omega |\nabla u|^p -\int_\Omega f u
&\leq {\ds \liminf_{\delta \rightarrow 0}}
\left(\frac{1}{p} \int_\Omega |\nabla u^\delta|^p- \int_\Omega f u^\delta\right)\\ 
&\leq {\ds \liminf_{\delta \rightarrow 0}} 
\left(\frac{1}{p} \int_\Omega |\nabla u^\delta|^p+ \frac{1}{\delta} 
\int_\Omega B(u^\delta-\psi) - \int_\Omega f u^\delta\right)\\ 
&\leq {\ds \liminf_{\delta \rightarrow 0}} \inf_{u \geq \psi} 
\left\{\frac{1}{p} \int_\Omega |\nabla u|^p+ \frac{1}{\delta} 
\int_\Omega B(u-\psi) - \int_\Omega f u \right\}\\
&\leq \frac{1}{p} \int_\Omega |\nabla v|^p -\int_\Omega f v.
\end{align*}
Then, one concludes that 
$ \nabla u^\delta \longrightarrow \nabla u \ {\rm{strongly \ in}} \ L^p(\Omega)$
and since $u \in K(\psi)$, then $u$ solves (\ref{P-ob:eq:1.1}).
 
\bigskip

(2)
(i) let $u^\delta$ be the solution of
(\ref{P-ob:eq:1.7}), since $\nabla u^\delta$ 
is uniformly bounded in $L^p(\Omega)$ by some constant $C$, we get that
$-\Delta_p u^\delta -f$ is bounded in $W^{-1,p^\prime}(\Omega)$, so it converges weakly, up to a subsequence, in $W^{-1,p^\prime}(\Omega)$.
Hence,  $-\frac{1}{\delta} \beta(u^\delta-\psi)$ converges too, up to a subsequence, in $W^{-1,p^\prime}(\Omega)$, and we have
$$-\frac{1}{\delta} \beta(u^\delta-\psi)  \rightharpoonup \mu \   
{\rm{weakly \ in}} \ W^{-1,p^\prime}(\Omega),$$
where $\mu$ is a positive distribution, hence a positive measure.
Then, by (1), we see that $u$ and $\mu$ are linked  by the relation (\ref{P-ob:eq:1.5}).

We now prove 
(ii): let $u$ be the solution of (\ref{P-ob:eq:1.1}). 
Taking $\phi=(\psi-u) \in W^{1,p}_0(\Omega) $ in the above inequalities,
one gets
$$- \frac{1}{\delta} \int_\Omega \beta(u^\delta-\psi) \ (u-\psi) \ dx 
\leq  \|\nabla u^\delta\|^{p-1}_p \|\nabla (\psi-u)\|_p
+\|f\|_{p^\prime} \| \psi-u\|_p .$$
Since $u \in K(\psi)$, passing to the limit we obtain:
$$\langle \mu, \psi-u \rangle = \int_\Omega |\nabla u|^{p-2} \nabla u
\cdot\nabla (\psi-u)- \int_\Omega f (\psi-u) = 0,\ \text{by} \ (\ref{P-ob:eq:1.2})$$
Then (ii) follows.
\end{proof}

\section{Optimal Control for a Non-Positive Source Term}
\subsection{Optimal control for a non positive source term}

\begin{proposition}\label{POC:prop:1.1}
Let $f, \psi$ and $T_f(\psi)$ be  as in {\rm{(\ref{P-ob:eq:1.1})}}. One has
$$\frac{1}{p}\int_\Omega |\nabla T_f(\psi)|^p \leq \frac{1}{p}\int_\Omega |\nabla \psi|^p \ dx+\int_\Omega f [T_f(\psi)- \psi] \  dx.$$
\end{proposition}

\begin{proof}[Proof of Proposition \ref{POC:prop:1.1}]
From {\rm{(\ref{P-ob:eq:1.1})}} taking $v=\psi$ and using H\"older's inequality, we have
$$\int_\Omega |\nabla T_f(\psi)|^p \leq 
\frac{p-1}{p} \|\nabla T_f(\psi)\|^p_p+
\frac{1}{p}\|\nabla \psi\|^p_p +\int_\Omega f [T_f(\psi)- \psi] \  dx.$$
\end{proof}

Note that since $T_f(\psi)\geq \psi$, it follows that if $f\leq 0$, then
\begin{equation}\label{POC:eq:1.1}
\int_\Omega |\nabla T_f(\psi)|^p \  dx \leq \int_\Omega |\nabla \psi|^p \  dx.
\end{equation}
Let us now introduce the following problem, said ``{\it{optimal control problem}}'':
\begin{equation}\label{POC:eq:1.2}
\inf_{\widetilde \psi \in W_0^{1,p}(\Omega)} J_f(\widetilde \psi),
\end{equation}
where
\begin{equation}\label{POC:eq:1.3}
J_f(\widetilde \psi)=\frac{1}{p}\int_\Omega \left\{|T_f(\widetilde \psi)- z|^p+|\nabla \widetilde \psi|^p \right\} \ dx,
\end{equation}
for some given $z \in L^p(\Omega)$. $z$ is said to be  {\it the
initial profile}, $\psi$ is {\it the control variable} and $T_f(\psi)$ is {\it the state variable}. The pair $(\psi^*,T_f(\psi^*))$ where $\psi^*$ is a solution
for (\ref{POC:eq:1.2}) is called an optimal pair and $\psi^*$ an optimal control.

In this section, we establish the existence and uniqueness of
the optimal pair in the case where $f\leq 0$. 

\begin{theorem}\label{POC:thm:1.1}
If $f \in L^{p'}(\Omega), \ f \leq 0$ on $\Omega$, then there exists
a unique optimal control $\psi^*\in W_0^{1,p}(\Omega)$ for 
{\rm{(\ref{POC:eq:1.2})}}.
Moreover, the corresponding state $u^*$ coincides with $\psi^*$,
i.e. $T_f(\psi^*)=\psi^*$.
\end{theorem}

\begin{proof}[Proof of Theorem \ref{POC:thm:1.1}]
In a first time we prove that there exists a pair of solutions of the form 
$(u^*, u^*)$, hence $(u^*=T_f(u^*))$.
Let $(\psi_k)_k$ be a minimizing sequence for (\ref{POC:eq:1.3}), then $T_f(\psi_k)$ is bounded 
in  $W^{1, p}(\Omega)$, therefore $T_f(\psi_k)$ converges for a subsequence towards some $u^* \in W_0^{1,p}(\Omega)$.
Moreover, using the lower semicontinuity of $T_f$ as in proposition \ref{P-ob:prop:1.5}, one gets
$$T_f(u^*) \leq \liminf_{k \rightarrow \infty}T_f(T_f(\psi_k))
\leq \lim_{k \rightarrow \infty}T_f(\psi_k)=u^*,$$
and by the definition of $T_f, \  T_f(u^*)\geq u^*$. 
Hence $u^*=T_f(u^*)$.

We prove that $(u^*,u^*)$ is an optimal pair. 
Using proposition \ref{P-ob:prop:1.5}, by the lower semicontinuity in $W_0^{1,p}(\Omega)$
of $T_f$:
\begin{align*}
J_f(u^*)
&=\frac{1}{p}\int_\Omega \left\{|u^*- z|^p+|\nabla u^*|^p \right\} \ dx\\
&\leq \liminf_{k \rightarrow \infty}
\frac{1}{p}\int_\Omega \left\{|T_f(\psi_k)- z)|^p+|\nabla \psi_k|^p \right\} \ dx\\
&=\inf_{\psi \in W_0^{1,p}(\Omega)} J_f(\psi).
\end{align*}
Secondly, we prove that every optimal pair is of the form $(u^*, u^*)$.
Observe that if $(\psi^*, T_f(\psi^*))$ is a solution then 
$(T_f(\psi^*), T_f(\psi^*))$ is a solution.
Indeed
$$\int_\Omega \left\{|T_f(\psi^*)- z)|^p+|\nabla T_f(\psi^*)|^p \right\} \ dx
\leq \int_\Omega \left\{|T_f(\psi^*)- z)|^p+|\nabla \psi^*|^p \right\} \ dx.$$
So
\begin{equation}\label{POC:eq:1.4}
\int_\Omega |\nabla T_f(\psi^*)|^p  \ dx= \int_\Omega |\nabla \psi^*|^p  \ dx,
\end{equation}
by inequality (\ref{POC:eq:1.1}),
using the H\"older's inequality, one obtains then
\begin{align*}
0
&\leq \int_\Omega f(\psi^*-T_f(\psi^*)) \ dx \\
&\leq \int_\Omega \sigma(T_f(\psi^*)) \cdot \nabla(\psi^*-T_f(\psi^*)) \ dx\\ 
&\leq \int_\Omega |\nabla T_f(\psi^*)|^{p-2}\nabla T_f(\psi^*)\cdot \nabla\psi^*-
\int_\Omega |\nabla T_f(\psi^*)|^p\\
&\leq \left(\int_\Omega |\nabla T_f(\psi^*)|^p\right)^{\frac{p-1}{p}}
\left(\int_\Omega |\nabla T_f(\psi^*)|^p\right)^{\frac{1}{p}}-
\int_\Omega |\nabla T_f(\psi^*)|^p=0,
\end{align*}
which implies
$$\int_\Omega |\nabla T_f(\psi^*)|^{p-2}\nabla T_f(\psi^*) \cdot \nabla\psi^*
-\int_\Omega |\nabla T_f(\psi^*)|^p=0.$$
Let us recall that by convexity, one has the following inequality
$$\frac{1}{p}\int_\Omega |\nabla \psi^*|^p+\frac{p-1}{p} \int_\Omega
|\nabla T_f(\psi^*)|^p- \int_\Omega |\nabla T_f(\psi^*)|^{p-2}\nabla
T_f(\psi^*) \cdot \nabla\psi^* \geq 0.$$
Then the equality holds and by the strict convexity, one gets $\nabla (\psi^*)=\nabla
(T_f(\psi^*)) \ \text{a.e.}$, hence $\psi^*=T_f(\psi^*)$.
Finally, we deduce from the two previous steps that the pair is unique.
Suppose that $(u_1, u_1)$ and $(u_2, u_2)$ are two solutions, and consider
$(\frac{u_1+u_2}{2}, T_f(\frac{u_1+u_2}{2}))$. We prove that it is
also a solution.
Indeed:
\begin{align*}
\int_\Omega \left|\frac{u_1+u_2}{2}- z\right|^p
&+\left|\nabla T_f(\frac{u_1+u_2}{2})\right|^p  \ dx\\
&\leq \int_\Omega \left|\frac{u_1+u_2}{2}- z\right|^p+\left|\nabla \left(\frac{u_1+u_2}{2}\right)\right|^p  \ dx\\
&\leq  \frac{1}{2} \left(J_f(u_1)+ J_f(u_2)\right)=
\inf_{\psi \in W_0^{1,p}(\Omega)} J_f(\psi),
\end{align*} 
which implies that $u_1=u_2$. Thus, the uniqueness of the optimal pair
for $f\leq0$ holds.
\end{proof}

\subsection{Optimal control for a nonnegative source term}
We are interested here to the case $f\geq 0$ on $\Omega$.
In what follows we will denote by $Gf$ the unique function in $W^{1, p}_0(\Omega)$ which verifies 
$$\begin{cases}
-\Delta_p (Gf)=f,\  \text{in} \ \Omega \ \text{a.e.}\\
Gf=0, \  \text{on} \ \partial \Omega,
\end{cases}$$
where $f \in L^{p'}(\Omega)$ and $Gf \in W^{1,p}_0(\Omega)$.

\begin{theorem}\label{POC:thm:1.2}
Suppose that $f\in L^{p^\prime}(\Omega)$ is a nonnegative function. Suppose that $z \in L^p(\Omega)$, satisfying $z\leq Gf$ a.e on $\Omega$. Then the minimizing problem {\rm{(\ref{POC:eq:1.2})}}
has a unique optimal pair $(0, Gf)$.
\end{theorem}

\begin{lemma}\label{POC:lemma:1.1}
Let $T_f(\psi)$ be a solution to {\rm{(\ref{P-ob:eq:1.1})}} and $Gf$ defined as above.
Then $T_f(\psi)$ is greater than $Gf$.
\end{lemma}

\begin{proof}[Proof of Lemma \ref{POC:lemma:1.1}]
We have that $-\Delta_p(Gf)=f$, and $T_f(\psi)$ realizes $-\Delta_p(T_f(\psi))\geq f$. 
Then, by the Comparison Theorem for $-\Delta_p$ we get that $Gf\leq T_f(\psi).$
\end{proof}

\begin{proof}[Proof of Theorem \ref{POC:thm:1.2}]
In a first time we prove that $(0,Gf)$ is an optimal pair. Indeed, for all $\psi \in W^{1,p}_0(\Omega)$
\begin{align*}
J_f(\psi)
&=\frac{1}{p} \int_\Omega \left\{|Gf-z+T_f(\psi)-Gf|^p+|\nabla \psi|^p\right\}\\
&\geq \frac{1}{p} \int_\Omega \left\{|Gf- z|^p +p|Gf- z|^{p-2}(Gf- z)(T_f(\psi)- Gf)\right\}\\
&\geq \frac{1}{p} \int_\Omega \left\{|Gf- z|^p\right\}\\
&=J_f(0).
\end{align*}
The equality with $(\psi^*, T_f(\psi^*))$ implies that we have equality in each step, so we get $\|\nabla \psi^*\|_p=0$, then $\psi^*=0 \ \text{a.e. \ in} \ \Omega$. 
Thus, $(0,Gf)$ is the unique optimal control pair. 
\end{proof}

\end{document}